\documentclass[preprint,12pt]{elsarticle}

\usepackage{amssymb}
\usepackage{amsmath}
\usepackage{amsthm}
\usepackage[colorlinks]{hyperref}
\AtBeginDocument{%
  \hypersetup{%
    linkcolor=blue,%
    citecolor=green,%
  }%
}
\usepackage{cleveref}
\usepackage{geometry}
\crefname{equation}{}{}

\newtheorem{theorem}{Theorem}[section]

\newtheorem{corollary}[theorem]{Corollary}

\theoremstyle{definition}

\theoremstyle{remark}

\numberwithin{equation}{section}

\journal{~}

\begin{document}

\begin{frontmatter}



\title{A simple proof of the transcendence of the trigonometric functions \tnoteref{t1}}

\author[rvt]{Yuanyuan Lian}
\ead{lianyuanyuan@nwpu.edu.cn}
\author[rvt]{Kai Zhang\corref{cor1}}
\ead{zhang\_kai@nwpu.edu.cn}
\tnotetext[t1]{Research is supported by NSFC 11701454.}

\cortext[cor1]{Corresponding author.}

\address[rvt]{Department of Applied Mathematics, Northwestern Polytechnical University, Xi'an, Shaanxi, 710129, PR China}

\begin{abstract}
In this note, we give a simple proof that the values of the trigonometric functions at any nonzero rational number are transcendental numbers.
\end{abstract}

\begin{keyword}
Irrational number \sep Transcendental number \sep Trigonometric function%  \sep H\"{o}lder continuity \sep Reifenberg flat domain \sep Fully nonlinear elliptic equation \sep Viscosity solution

\MSC[2010] 11J81

\end{keyword}

\end{frontmatter}


%\section{Introduction}
%\label{S1}

The irrationality and transcendence of numbers have attracted a lot of interest, especially for the numbers related to $e$ and $\pi$. In 1737, Euler \cite{Euler} proved the irrationality of $e$. The irrationality of $\pi$ was proved by Lambert \cite{Lambert} (see also \cite[Chap. 34]{MR0121327}) in 1761. Both proofs depend heavily on the continued fractions. In 1873, Hermite \cite{Hermite} gave the proof of the transcendence of $e$ by a new technique, which involves constructing an auxiliary function and the integral by parts. This technique has been developed extensively and all the results cited below are based on this technique. The transcendence of $\pi$ was finally proved by Lindemann \cite{Lindeman1,Lindeman2} in 1882.

In addition, it can be proved (see \cite{MR0028344, MR0027784} and \cite[Theorem 48, Chap. 3]{MR2445243}) that the exponential function $e^x$ maps rational numbers ($\neq 0$) to irrational numbers. In 1947, Niven \cite{MR0021013} gave a simple proof of the irrationality of $\pi$. A modified proof shows that the trigonometric functions map rational numbers ($\neq 0$) to irrational numbers (see \cite[Theorem 2.5]{MR0080123}). With the aid of the idea of Hurwitz \cite{MR1510809}, Niven \cite{MR0080123} gave a simple proof of the transcendence of $e$. It indicates that the exponential function $e^x$ maps rational numbers ($\neq 0$) to transcendental numbers.

In this note, based on the technique of Hermite, we give a simple proof of the transcendence of trigonometric functions. That is, the trigonometric functions map rational numbers ($\neq 0$) to transcendental numbers.

\begin{theorem}\label{t-1}
For any rational number $r\neq 0$, $\cos r$ is a transcendental number.
\end{theorem}

\noindent \textbf{Proof.} Since $\cos r=\cos (-r)$, we can assume that $r$ is positive and write $r=s/t$ where $s,t$ are relatively prime positive integers.

We use a proof by contradiction. Suppose that $\cos r$ is not a transcendental number. Then there exist a positive integer $m$ and rational numbers $\bar{a}_j$ ($0\leq j\leq m$) such that
\begin{equation*}
  \sum_{j=0}^{m}\bar{a}_j\cos^j r=0.
\end{equation*}
From the product-to-sum formulas for trigonometric functions, there exist rational numbers $a_j$ ($0\leq j\leq m$) such that
\begin{equation}\label{e.1}
\sum_{j=0}^{m} a_j \cos jr=0.
\end{equation}
By multiplying an integer, we may assume that $a_j$ ($0\leq j\leq m$) are all integers.

If $a_0 \neq 0$, let
\begin{equation}\label{fx}
  f(x)=\frac{t^{4mp+2p-2}x^{2p-2}(x^2-r^2)^{2p}(x^2-4r^2)^{2p} \cdots (x^2-m^2r^2)^{2p}}{(2p-2)!},
\end{equation}
where $p$ is a positive prime to be specified later. If $a_0=0$ and $a_{j_0}\neq 0$ for some $1\leq j_0\leq m$, we take $f$ to be
\begin{equation*}
\begin{aligned}
f(x)=\frac{t^{4mp+2p-2}f_0(x)f_1(x)\cdots f_{j_0-1}(x)(x-j_0r)^{2p-2}f_{j_0+1}(x)\cdots f_m(x)}{(2p-2)!}
\end{aligned}
\end{equation*}
where $f_k(x)=\left((x-j_0r)^2-(kr-j_0r)^2\right)^{2p}$ for $k\neq j_0$. In both cases, the following arguments are almost the same. Thus, we only give the proof for the case $a_0\neq 0$. For $0<x<mr$, we have
\begin{equation}\label{abfx}
  |f(x)|<\frac{t^{4mp+2p-2}(mr)^{4mp+2p-2}}{(2p-2)!}.
\end{equation}
It can also be easily verified that for any $k\geq 0$, $f^{k}(0)$ is an integer and $f^{(2k+1)}(0)=0$ (see \cite[Chapter 2.2]{MR0080123}).

By the elementary calculus, for any positive integer $j$,
\begin{equation*}
  \int_{0}^{jr} f(x)\sin(jr-x)dx=
  \sum_{k=0}^{\infty}(-1)^{k}f^{(2k)}(jr)+
  \sum_{k=0}^{\infty}(-1)^{k}f^{(2k)}(0)\cos jr+
  \sum_{k=0}^{\infty}(-1)^{k}f^{(2k+1)}(0)\sin jr.
\end{equation*}
Note that $f^{(2k+1)}(0)=0$ for any $k$. Hence,
\begin{equation*}
    \int_{0}^{jr} f(x)\sin(jr-x)dx=
  \sum_{k=0}^{\infty}(-1)^{k}f^{(2k)}(jr)+
  \sum_{k=0}^{\infty}(-1)^{k}f^{(2k)}(0)\cos jr.
\end{equation*}
Then, by\cref{e.1},
\begin{equation}\label{msum}
  \sum_{j=0}^{m}a_{j}\int_{0}^{jr} f(x)\sin(jr-x)dx=\sum_{j=0}^{m}\sum_{k=0}^{\infty}(-1)^{k}a_{j}f^{(2k)}(jr).
\end{equation}

The right hand in above equation is an integer. Moreover, by\cref{fx},
$f^{(2k)}(jr)$ is divisible by $p$ for all $k$ and $j$ with one exception:
\begin{equation*}
  f^{(2p-2)}(0)=s^{4mp}t^{2p-2}(m!)^{4p},
\end{equation*}
if we choose $p>\max\left\{m,s,t\right\}$. Next, by taking $p>|a_0|$, the right hand of\cref{msum} consists of a sum of multiples of $p$ with one exception, namely $a_0 f^{(2p-2)}(0)$. Hence, the right hand of\cref{msum} is a non-zero integer. However, by\Cref{abfx}, the left hand of\Cref{msum} satisfies
\begin{equation*}
  \left|\sum_{j=0}^{m}a_j\int_{0}^{jr}f(x)\sin(jr-x)dx \right|
  \leq \sum_{j=0}^{m}|a_j|\cdot jr \cdot \frac{t^{4mp+2p-2}(m^2r^2)^{2mp+2p-2}}{(2p-2)!}
  <1,
\end{equation*}
provided $p$ is chosen sufficiently large. Thus we have a contradiction, and the theorem is proved.~\qed\\

\begin{corollary}
  The trigonometric functions are transcendental at non-zero rational values of the arguments.
\end{corollary}
\noindent \textbf{Proof.} Because of $\cos 2r=1-\sin^2 r$ and $\cos 2r=(1-\tan^2r)/(1+\tan^2r)$, $\sin r$ and $\tan r$ are transcendental numbers for any rational number $r\neq 0$. Also, $\csc r$, $\sec r$ and $\cot r$ are transcendental numbers for any rational number $r\neq 0$.~\qed

\section*{References}
\bibliographystyle{model2-names}
\bibliography{Cosr}

\end{document}